\documentclass{article}%

\usepackage{amsfonts}
\usepackage{amsmath}
\usepackage{amsthm,amstext,amsfonts,bm,amssymb,amsopn}
\usepackage{amssymb}
\usepackage{xfrac}
\usepackage{graphicx}%
\usepackage{color}%
\usepackage{refcount}
\allowdisplaybreaks
\usepackage{xypic}
\usepackage{environ,enumitem}
\usepackage[margin=.5cm]{caption}

\numberwithin{equation}{section}

\makeatletter

\usepackage{tikz}

\makeatletter
\newcommand{\xleftrightarrow}[2][]{\ext@arrow 3359\leftrightarrowfill@{#1}{#2}}
\makeatother

\newcommand{\xdasharrow}[2][->]{
\tikz[baseline=-\the\dimexpr\fontdimen22\textfont2\relax]{
\node[anchor=south,font=\scriptsize, inner ysep=1.5pt,outer xsep=5pt](x){#2};
\draw[shorten <=3.4pt,shorten >=3.4pt,dashed,#1](x.south west)--(x.south east);
}
}

\providecommand{\U}[1]{\protect\rule{.1in}{.1in}}
\providecommand{\comment}[1]{}

\makeatletter
\newtheorem*{rep@theorem}{\rep@title}
\newcommand{\newreptheorem}[2]{%
\newenvironment{rep#1}[1]{%
 \def\rep@title{#2 \ref{##1}}%
 \begin{rep@theorem}}%
 {\end{rep@theorem}}}
\makeatother

\newtheorem*{namedtheorem}{\theoremname}
\newcommand{\theoremname}{testing}

\newtheorem{theorem}{Theorem}
\newtheorem*{theorem*}{Theorem}
\newtheorem*{fact*}{Fact}

\newtheorem{proposition}[theorem]{Proposition}

\theoremstyle{definition}

\newreptheorem{theorem}{Theorem}
\newreptheorem{lemma}{Lemma}
\newreptheorem{proposition}{Proposition}

\NewEnviron{thmquote}{\vspace{1ex}\par
\refstepcounter{equation}%
\hfill\parbox{\dimexpr \textwidth-2cm}%
{\centering\small\textit{\BODY}}%
\hfill\llap{(\thequote)}\vspace{1ex}\par}

\newcommand{\BR}{\mathbb R}

\title{Metrizing the Chabauty topology}
\author{Ian Biringer}

\begin{document}

\maketitle

\begin {abstract}
 We describe an explicit metric that induces the Chabauty topology  on the space of closed subsets of a proper metric space $M$.
\end {abstract}

\section{Introduction}

Suppose $M $ is a proper metric space and let $\mathcal C (M) $ be the space of closed subsets of $M $.  The \emph {Chabauty topology} on $\mathcal C(M)$ is generated by the subsets 
\begin {equation}\label {neighborhoods}\{C \in \mathcal C (M) \ | \ C \cap K=\emptyset\}, \ \   \  \ \{C \in \mathcal C (M) \ | \ C \cap U \neq\emptyset\},
\end {equation}
where $K\subset M$ is compact and $U\subset M $ is open.  See Chabauty \cite{Chabautylimite} and \cite[Ch E]{Benedettilectures}.

When $M$ is compact, the Chabauty topology is induced by the \emph {Hausdorff metric} on $C(M)$, where the distance between closed subsets $C_1,C_2 \subset M$ is $$d_{\text {Haus}}(C_1,C_2)=\inf \{ \epsilon \ | \ C_1 \subset \mathcal N_\epsilon (C_2) \text { and }  C_2 \subset \mathcal N_\epsilon (C_1) \}. $$
In general, it is well-known that the Chabauty topology is compact, separable and metrizable \cite[Lemma E.1.1]{Benedettilectures}, but  most of the metrizability proofs in the literature go through Urysohn's theorem. 

 In search of an explicit metric, note that the Chabauty topology is almost, but not quite, induced by taking the Hausdorff topology on all compact subsets of $M $.  Namely, fix a base point $p\in M$.  If $A\subset M$ is closed and $R>0$, set
$$
A_R=A\cap \overline {B(p,R)},$$
and then define a pseudo-metric $d_R$ on $\mathcal C (M) $ by setting
$$ d_R(A,B)=\min  \Big \{ 1, d_{\text {Haus}}(A_R, B_R) \Big\},$$
where $d_{\text {Haus}}$ is the Hausdorff metric of the compact subset $\overline {B(p,R)} \subset M $. 

The family of pseudo-metrics $\{d_R \ | \ R>0\}$ does not determine the Chabauty topology, since if $x_i\rightarrow x$ is a convergent sequence of points with $d(x,p)=R$ and $d(x_i,p)>R$ for all $i $, then $\{x_i\} \rightarrow \{x\}$ in the Chabauty topology, but $d_R(\{x_i\},\{x\})=1$ for all $i$.  However, the following is true:

\begin {theorem}\label {Chabautymetric}
The Chabauty topology on $\mathcal C(M)$ is induced by the metric
$$d : \mathcal C(M) \times \mathcal C(M) \longrightarrow\BR ,  \ \ d(A,B)=\int_0^\infty e^{-R} d_R(A,B) \, dR.$$
\end {theorem}

The point is that a Chabauty convergent sequence can fail to $d_R$-converge for only countably many $R $, a discrepancy which disappears under integration.  Above, $e^{-R}$  could be replaced by any positive, integrable function on $[0,\infty).$

\vspace{2mm}

We imagine that those who are sufficiently interested could probably come up with more metrics inducing the Chabauty topology (with Ab\'ert, we produce  a different one in \cite[A.4]{Abert_Biringer}).    However, we have not seen the expression above in the literature, and we think that  the way it formalizes the intuition that the Chabauty topology is almost the ``Hausdorff topology on compact sets'' is beautiful enough to justify this short note.


\subsection{Acknowledgements}

The author  is partially supported by NSF grant DMS 1611851. Thanks are due to the referee for improving the clarity of the paper.

\section{Proof of Theorem \ref {Chabautymetric}}

Before beginning the proof, recall that convergence in the Chabauty topology  can be characterized as follows.

\begin {proposition}[Prop E.12, \cite{Benedettilectures}]\label {Chabautyconvergence}
A sequence $(C_i)$ in $\mathcal C (M) $ converges to $C \in \mathcal C (M)$ in the Chabauty topology if and only if
\begin {enumerate}
\item if $x_{i_j} \in C_{i_j}$ and $x_{i_j} \to x \in M$, where $i_j \to \infty$, then $x\in C$.
\item if $x\in C$, then there exist $x_i \in C_i$ such that $x_i \to x$.
\end {enumerate}
\end {proposition}

Let $d$ be  the metric in Theorem \ref {Chabautymetric}. As the Chabauty topology is first countable, it suffices to show that a sequence Chabauty-converges if and only if it $d$-converges.

 Suppose that $(C_i)$ converges to $C$ with respect to $d$. Then $d_R(C_i,C)\rightarrow 0 $ for a.e.\ $R$, so in particular for arbitrarily large $R $. We check that $(C_i)$ Chabauty converges to $C $ using Proposition \ref {Chabautyconvergence}.  With $x$ as defined therein, we can just take any $R>d(x,p)$ with $d_R(C_i,C)\rightarrow 0$ and use that $d_R$ defines the Chabauty topology on the ball $B(p,R)$ to say that 1) and 2) are satisfied.

On the other hand, suppose that $(C_i)$ Chabauty-converges to $C $. We claim that $d_R (C_i, C)\rightarrow 0 $ for all but countably many $R $.  This will suffice to prove the proposition, for as each $d_R \leq 1$, the functions $e^{-R}d_R(C_i,C)$ are bounded by the integrable function $e^{-R}$, so then we must have $$d(C_i,C)=\int_0^\infty e^{-R} d_R(C_i,C) \, dR \rightarrow 0$$ by the Dominated Convergence Theorem.

So, our goal is to show that for all but countably many $R $, we have \begin {equation}\label {convsub} C_i\cap  \overline {B(p,R)} \rightarrow C\cap \overline {B(p,R)}\end {equation}
in the Chabauty topology on subsets of the compact set $ \overline {B(p,R)}$. In light of Proposition~\ref{Chabautyconvergence},  we need to show:
\begin {enumerate}
\item[1'.] if $x_{i_j} \in C_{i_j} \cap \overline {B(p,R)}$ and $x_{i_j} \to x \in \overline {B(p,R)}$, where $i_j \to \infty$, then $x\in C \cap \overline {B(p,R)}$.
\item[2'.] if $x\in C\cap \overline {B(p,R)}$, then there exist $x_i \in C_i \cap \overline {B(p,R)}$ such that $x_i \to x$.
\end {enumerate}

 It follows from the Chabauty convergence $C_i \rightarrow C$ that property 1' holds for every $R$. So, the point is to  prove 2' for all but  countably many $R$.

We claim that 2' holds when $R$ is chosen so that every $x\in C$ with $d(x,p)=R$ is in the closure of $C\cap B(p,R)$. For if $$x\in C\cap \overline {B(p,R)},$$ either $x\in B(p,R)$ or $d(p,x)=R$. In the first case, the Chabauty convergence $C_i \rightarrow C$ implies that there is a sequence $x_i \in C_i$ converging to $x$; eventually, these $x_i \in B(p,R)$, so we're done. In the second case, we know $x$ is the limit of a sequence $y_n \in C\cap B(p,R)$. Each $y_n$ is the limit of a sequence $y_{n,i} \in C_i\cap  \overline {B(p,R)}$, as in the first case, and then $y_{i,i}\rightarrow x$. So, 2' holds.

Finally, we claim that the condition in the previous paragraph fails for only countably many $R $, which amounts to proving that if
$$L=\Big\{ x\in C \ \Big| \ x\notin \overline {C\cap B(p,d(p,x))} \Big\},$$
then the set $\{d(x,p) \ | \ x \in L\}$ is countable.
For every point $x\in L$, there is some $\epsilon(x) >0 $ such that $$C \cap B(p,d(p,x))\cap B(x,\epsilon(x))=\emptyset.$$ We claim that for every compact $K \subset M$ and $\epsilon_0>0$, there are only finitely many values $d(p,x)$ where $x\in L\cap K$ and $\epsilon (x)\geq\epsilon_0 $. Exhausting $M $ with a countable union of compact sets and taking a (countable) sequence of such $\epsilon_0$ converging to $0$ will prove $L$ is countable.

If there are not finitely many such values $d(p,x)$, there is an infinite sequence of points $x_i \in L\cap K$ with $d(p,x_i)$ all distinct, and $\epsilon (x_i)\geq \epsilon_0 $. Passing to a subsequence, we may assume that $x_i \rightarrow x\in K$. Then for large $i,j$, we have
\begin{equation}
	x_i \in C \cap B(x_j,\epsilon_0) \ \text{ and } \ x_j \in C \cap B(x_i,\epsilon_0).\label{lasteq}
\end{equation}
Moreover, since the  distances $d(p,x_i)$ and $d(p,x_j)$ are distinct, we have either
\begin{equation} 
	x_i \in B(p,d(p,x_j)) \  \text{ or } \ x_j \in B(p,d(p,x_i)),\label{lasteq2}
\end{equation}
 so combining \eqref{lasteq} and \eqref{lasteq2} we contradict the definition of either $\epsilon(x_i)$ or $\epsilon(x_j)$, since both are at least $\epsilon_0$. 

\textit{\textrm{
\bibliographystyle{amsplain}
\bibliography{bibrefs}
}}

\end{document}